\documentclass[a4paper,11pt]{article}

\usepackage[utf8]{inputenc}
 \usepackage[T1]{fontenc}
 \usepackage[normalem]{ulem}
 \usepackage[english]{babel}
\usepackage{amsmath}
  \usepackage{amsthm}
 \usepackage{bbm}
 \usepackage{amssymb}
 \usepackage{verbatim}
 \usepackage{vmargin}
 \usepackage{graphicx}
 \usepackage{color}
 \usepackage{epsfig}
 \usepackage[only,llbracket,rrbracket]{stmaryrd}

\usepackage{setspace}

\newtheorem{thm}{Theorem}[section]
\newtheorem*{deft}{Definition}
\newtheorem*{notat}{Notations}
\newtheorem{prop}{Proposition}[section]
\newtheorem{lem}{Lemma}[section]
\newtheorem{coro}{Corollary}[section]
\numberwithin{equation}{section}

\theoremstyle{remark}
\newtheorem*{rques}{\textbf{Remarks}}{\vskip 0.5cm}

\newtheorem*{rque}{\textbf{Remark}}{\vskip 0.5cm} 
\newtheorem*{ack}{\textbf{Acknowledements}}{\vskip 0.5cm} 

\title{The three-dimensional \\ finite Larmor radius approximation}
\author{Daniel Han-Kwan\\Département de Mathématiques et Applications \\ Ecole Normale Supérieure \\ 45 rue d'Ulm \\ 75230 Paris Cedex 05, France \\Tel. No.: +33 1 44 32 20 65 \\Fax No.: +33 1 44 32 20 80\\E-mail: hankwan@dma.ens.fr}

\begin{document}
 \maketitle

\begin{abstract}
 Following Frénod and Sonnendrücker (\cite{FS2}), we consider the finite Larmor radius regime for a plasma submitted to a large magnetic field and take into account both the quasineutrality and the local thermodynamic equilibrium of the electrons. We then rigorously establish the asymptotic gyrokinetic limit of the rescaled and modified Vlasov-Poisson system in a three-dimensional setting with the help of an averaging lemma.
\end{abstract}
 

 \section{Introduction and main results}
 
 \subsection{Physical motivation}
 
 We are interested in the behaviour of a plasma (\textit{id est} a gaz made of ions with individual charge $Ze$ and mass $m_i$ and electrons with individual charge $-e$ and mass $m_e$, with $m_i >> m_e$) which is submitted to a large external magnetic field. It is ``well-known'' that such a field induces fast small oscillations for the particles and consequently introduces a new small time scale which is very restrictive and inconvenient from the numerical point of view. The simulation of such plasmas appears to be primordial since the model can be applied to tokamak plasmas from magnetic confinement fusion (like for the ITER project).

 \subsubsection{Heuristic study}
 Let us give some heuristic formal arguments to investigate the behaviour of the plasma:
 if we consider the motion of one particle (of charge $q>0$, mass $m$, position $x$ and velocity $v$) submitted to an external constant field $B$, the fundamental principle of mechanics gives that:
 \begin{equation}
 \frac{dx}{dt}=v, \text{                          }\frac{dv}{dt}=\frac{q}{m}( v\wedge  B)
 \end{equation}
 Straightforward calculations show first of all that the parallel velocity, denoted by $v_\parallel$ (that is to say the component of the velocity in the direction of the magnetic field) is conserved and thanks to the conservation of the kinetic energy, so is the norm of the perpendicular velocity $v_\perp$ (the component of the velocity in the perpendicular plane). 
 Actually, we can see that the particle moves on a helix whose axis is the direction of the magnetic field.
 The rotation period (around the axis) is the inverse of the cyclotron frequency $\Omega$:
 
  \begin{equation}
\Omega=\frac{\vert q\vert \vert B\vert  }{m}
  \end{equation}
  and the radius is the so-called Larmor radius:
 \begin{equation}
 r_L= \frac{\vert v_\perp \vert}{\Omega}
  \end{equation}
  
  In the case where the magnetic field is very strong, $\Omega$ tends to infinity whereas $r_L$ tends to zero. More precisely, if we take $\vert B \vert \sim \frac{1}{\epsilon}$ (with $\epsilon \rightarrow 0$) we have:

$$
\left\{
    \begin{array}{ll}
  \Omega \sim \frac{1}{\epsilon} \\
  r_L \sim \epsilon
\end{array}
  \right.
$$
  
  The approximation which consists in considering $r_L=0$ is the classical guiding center approximation (\cite{Grad}). This means that each particle is assimilated to its ``guiding center'' (in other words its ``instantaneous rotation center''), which is equivalent to neglect the very fast rotation of the particle around the axis.
  
  If one also applies some external constant electric field $E$, a similar computation shows that there appears:
  \begin{enumerate}
  \item an acceleration  $\frac{E.B}{\vert B\vert}$ in the direction of $B$. If we consider $E\sim 1$, then: 
\begin{equation}\frac{E.B}{\vert B\vert}\sim 1\end{equation}
  \item a drift $\frac{E\wedge B}{\vert B\vert^2}$ in the orthogonal plane. We have:
  
  \begin{equation}
  \frac{E\wedge B}{\vert B\vert^2} \sim \epsilon
  \end{equation}
  
  \end{enumerate}
  This drift, usually called the electric drift is problematic as regards to the issue of plasma confinement. It is negligible compared to the acceleration in the direction of $B$, but in the time scale for plasma fusion which is expected to be very long, one can not neglect this small drift, since it creates a displacement of order $\epsilon t$ ($t$ represents the time).

  At last, note also that if the fields are not constant, various other drifts may appear, whose order in $\epsilon$ is higher than those of the electric drift.

  Actually, the fields considered are neither constant, nor external, but self-induced by the plasma itself. The effects we would like to describe are due to the non-linear interaction between the particles and the electromagnetic field. 

\subsubsection{The mathematical model}

In all the sequel, we assume that the magnetic field is external and constant and we suppose that the speed of particles is small compared to the speed of light, so that we can use the electrostatic approximation which consists in reducing the Maxwell equations to the Poisson equation.  Finally, we decide to opt for a kinetic description for the ions: in other words, the time and space scales considered here are such that ions are not at a thermodynamic equilibrium and their density is governed by a kinetic equation.

  The basic model usually considered for the ions is the following Vlasov-Poisson system:
$$
\left\{
    \begin{array}{ll}
  \partial_t f + v.\nabla_x f + (E+ v\wedge B).\nabla_v f = 0 \\
  E= -\nabla_x V \\
  -\Delta_x V = \int f dv \\
f_{t=0}=f_0
\end{array}
  \right.
$$

  where $f(t,x,v)$ is the density of ions, with $t\in \mathbb{R}^+, x \in \mathbb{R}^d \text{ or } \mathbb{R}^d/\mathbb{Z}^d , v \in \mathbb{R}^d$ (usually $d=2$ or $3$), meaning that $f(t,x,v)dxdv$ gives the number of ions in the infinitesimal volume $[x,x+dx]\times [v,v+dv]$ at time $t$ (note that in this model, electrons are for the moment neglected).

  \subsubsection{The gyrokinetic approximation}
  
   It is important from a numerical point of view to establish the asymptotic equation when $\vert B \vert$ tends to infinity. Indeed, we expect the asymptotic equation to be ``easier'' to handle: only one time and space scale, perhaps less variables in the phase space to deal with... The derivation of such equations is usually referred to in the mathematic literature as ``gyrokinetic approximation''.
   
  Rigorous justifications of these derivations with various time and space observation scales have only appeared at the end of the nineties. We refer for instance to the works of Brenier (\cite{Br}), Frénod and Sonnendrücker (\cite{FS1}-\cite{FS2}), Frénod, Raviart and Sonnendrücker (\cite{FRS}), Golse and Saint-Raymond (\cite{GSR1}-\cite{GSR2}), Saint-Raymond (\cite{SR1}-\cite{SR2}).
  
  The classical ``guiding center approximation'' corresponds to the following scaling for the Vlasov-Poisson system (from now on and until the end of the paper, $B$ is a constant vector, say for instance $B=\frac{1}{\epsilon}e_z$):
  \begin{equation}
   \label{VPbase}
\left\{
    \begin{array}{ll}
    \partial_t f_\epsilon + v.\nabla_x f_\epsilon + (E_\epsilon+ \frac{v\wedge e_z}{\epsilon}).\nabla_v f_\epsilon = 0 \\
  E_\epsilon= -\nabla_x V_\epsilon \\
  -\Delta_x V_\epsilon = \int f_\epsilon dv \\
     f_{\epsilon, t=0}=f_0
\end{array}
  \right.
\end{equation}

  The articles \cite{FS1} and \cite{GSR1} show that when $\epsilon \rightarrow 0$, this leads to a one-dimensional kinetic equation in the direction of $B$:

 \begin{equation}
\left\{
    \begin{array}{ll}  
\partial_t f + v_\parallel.\nabla_x f + E_\parallel.\nabla_v f = 0 \\
  E= -\nabla_x V \\
  -\Delta_x V = \int f dv \\
     f_{t=0}=f_0 \\
\end{array}
  \right.
\end{equation}
Notice that the electric drift does not appear; this was expected since we have seen in the formal analysis that this drift was of higher order in $\epsilon$ than the other effects. This shows in particular that this approximation is not sufficient for the numerical simulation of tokamaks. In order to make this drift appear, there exists to our knowledge two main possibilities:
\begin{enumerate}
 \item  one consists in restricting to a  $2D$ problem in the plane orthogonal to $B$ (\cite{GSR1}), 
 \item the other consists in rescaling the orthogonal scales in order to get both transport and electric drift at the same order (\cite{FS2}).
\end{enumerate}

  This work directly follows the articles \cite{FS2} and \cite{FRS} where the authors considered the ``finite Larmor radius approximation''. This means that the spatial observation scale in the plane orthogonal to $B$ is chosen smaller than the one in the parallel direction, more precisely with the same order as the Larmor radius $r_L$, so that one can expect the electric drift to appear in the asymptotic equation. 

In some sense, having such a scaling allows the electric field to significantly vary across a Larmor radius, which is not the case for instance in (\ref{VPbase}).
Moreover, in this situation, the positions of the particles are no longer assimilated to the position of their ``guiding center'' and we will have to perform an average over one fast oscillation period (the so-called gyroaverage) in order to get a sort of averaged number density.

 \subsection{Scaling and existing results}

The system we are going to study is based on the ``finite Larmor radius scaling'' and takes into account the quasineutrality of the plasma.
\subsubsection{The (refined) mathematical model}

 We refer to \cite{FS2} for a complete discussion on the scaling. Let us recall briefly and quite crudely how it works. 

Let $L_\parallel$ be the characteristic length in the direction of the magnetic field and $L_\perp$ be the characteristic length in the perpendicular plane. We consider that $L_\parallel\sim 1$ and $L_\perp\sim \epsilon$ and define the dimensionless variables $x'_\parallel= \frac{x_\parallel}{L_\parallel}$ and $x'_\perp= \frac{x_\perp}{L_\perp}$. In the same fashion we also define the dimensionless variables $t'$ and $v'$ with characteristic time and velocities with the same order as $L_\parallel$ and introduce the new number density $f'$ defined by $\bar{f}f'(t',x',v')=f(t,x,v)$ (and we define likewise the new electric field and potential $\bar{E}E'(t',x',v')=E(t,x,v)$ and $\bar{V}V'(t',x',v')=V(t,x,v)$). 
We consider the scaling $\bar{f},\bar{E}\sim 1$ and $\bar{V}\sim {\epsilon}$.
At last, we introduce the Debye length of the plasma $\lambda_D$, which appears in the Poisson equation. In order to take into account the quasineutrality of the plasma, we take from now on  $\lambda_D \sim \sqrt{\epsilon}$. 

The Poisson equation states in this scaling:
\begin{eqnarray}
-\epsilon \Delta_{x_\parallel'} V'_\epsilon -  \frac{1}{\epsilon}\Delta_{x_\perp'} V_\epsilon' = \frac{1}{\epsilon}\left( n^i_\epsilon - n^e_\epsilon\right) 
 \end{eqnarray}
where $n^i_\epsilon = \int f'_\epsilon dv'$ is the density of ions and $n^e_\epsilon$ the density of electrons. 
The density distribution of ions is normalized so that $\int f'_0 dv' dx' =1$.

The main difference between Frénod and Sonnendrücker's model and ours lies in the following. Instead of considering a fixed background of electrons, and since $\frac{m_e}{m_i}<<1$, we make the usual assumption that the (adiabatic) electrons are instantaneously at a local thermodynamic equilibrium, so that their density follows a Boltzmann-Maxwell distribution:

\begin{equation}
 n^e_{\epsilon}(x,t)=\exp \left(\frac{e {V}'_\epsilon}{k_B T_e}\right) 
\end{equation}
where $k_B$ is the Boltzmann constant, $-e$ the charge and $T_e$ the temperature of the electrons. We consider that $\frac{e}{k_B T_e}\sim 1$.

We make the assumption that we are not far from a fixed background of electrons, so that we can linearize this expression:
\begin{equation}
 n^e_\epsilon(x,t)=1+ V'_\epsilon
\end{equation}
 
We are obviously aware that this assumption is not really satisfactory both from a mathematical and physical point of view; let us just consider this as a plain mathematical model which allows us to give rigorous justifications in this work. The problem of a fixed background of electrons, i.e. $n^e_\epsilon=1$, brings actually more interesting formal results; this point will be discussed in the last section.

The Poisson equation can now be written:
\begin{equation}
V'_\epsilon- \epsilon^2 \Delta_{x_\parallel'} V'_\epsilon - \Delta_{x_\perp'} V'_\epsilon = \int f_\epsilon dv' - \int f_0 dv'dx'
\end{equation}

The dimensionless system (\ref{VPbase}) becomes (for the sake of simplicity, we forget the primes):
  \begin{equation}
\label{real3D}
\left\{
    \begin{array}{ll}
    \partial_{t} f_\epsilon + \frac{v_\perp}{\epsilon}.\nabla_{x} f_\epsilon + v_\parallel.\nabla_{x} f_\epsilon + (E_\epsilon+ \frac{v\wedge B}{\epsilon}).\nabla_{v} f_\epsilon = 0 \\
  E_\epsilon= (-\nabla_{x_\perp} {V}_\epsilon, -\epsilon\nabla_{x_\parallel} {V}_\epsilon)  \\
  V_\epsilon- \epsilon^2 \Delta_{x_\parallel} V_\epsilon - \Delta_{x_\perp} V_\epsilon = \int f_\epsilon dv - \int f_0 dvdx \\
     f_{\epsilon, t=0}=f_{\epsilon, 0}
\end{array}
  \right.
\end{equation}
with the notation $\Delta_{x_\parallel}= \partial^2_{x_\parallel}$ and $\Delta_{x_\perp}=\Delta-\Delta_{x_\parallel}$,

the problem being posed for $(x_\perp,x_\parallel,v) \in \mathbb{T}^2\times \mathbb{T} \times \mathbb{R}^3$ (with $\mathbb{T}=\mathbb{R} / \mathbb{Z}$ equipped with the restriction of the Lebesgue measure to $[0,1[$).

\subsubsection{State of the art about the Finite Larmor Radius Approximation}
 
Using homogenization arguments, Frénod and Sonnendrücker established the convergence in some weak sense of sequences of solutions $(f_\epsilon)_{\epsilon \geq0}$ of similar systems, in two cases, namely in some pseudo $2D$ case (assuming that nothing depends on $x_\parallel$ and $v_\parallel$) and in a $3D$ case when the electric field is external. The main tool used to establish the convergence is the ``2-scale convergence'' introduced by Nguetseng \cite{Ngu} and Allaire \cite{Al} that we will recall later on.

\begin{enumerate}
 \item \underline{The 3D case}:

Assume that we deal with an external electric field $E_\epsilon = E \in \mathcal{C}^1(\mathbb{R}\times\mathbb{R}^3)$: 
$$
\left\{
    \begin{array}{ll}
   \partial_t f_\epsilon + v_\parallel .\nabla_x f_\epsilon +  \frac{v_\perp}{\epsilon} .\nabla_x f_\epsilon + \left( E + \frac{v\wedge e_z}{\epsilon} \right).\nabla_v f_\epsilon =0 \\
f_{t=0}=f_0
\end{array}
  \right.
$$
Frénod and Sonnendrücker proved the following theorem:
 \begin{thm} For each $\epsilon$, let $f_\epsilon$ be the unique solution of the scaled Vlasov equation in $L^\infty_t(L^1_{x,v}\cap L^2_{x,v})$. Then the following convergence holds as $\epsilon$ tends to $0$:
 \begin{equation}
 f_\epsilon  \rightharpoonup f  \text{    weak-*    } L^\infty_t(L^2_{x,v})
 \end{equation}
 where $f \in L^\infty_t(L^2_{x,v})$ is the unique solution to:
 \begin{eqnarray*}
  \partial_t f + v_\parallel .\nabla_x f+ \frac{1}{2\pi}\left(\int_0^{2\pi} \mathcal{R}(\tau)E(t,x+\mathcal{R}(-\tau)v)d\tau \right).\nabla_x f \\
  +\frac{1}{2\pi}\left(\int_0^{2\pi} R(\tau)E(t,x+\mathcal{R}(-\tau)v)d\tau \right).\nabla_v f =0
 \end{eqnarray*}
 $$  f_{\vert t=0} =   \frac{1}{2\pi}\left(\int_0^{2\pi} f_0(x+\mathcal{R}(\tau)v, R(\tau)v)d\tau \right)$$
 denoting by $R$ and $\mathcal{R}$ the linear operators defined by:
 $$R(\tau)= \begin{bmatrix}\cos \tau & -\sin \tau & 0 \\ \sin \tau & \cos \tau & 0 \\ 0 & 0 & 1\end{bmatrix}, \mathcal{R}(\tau) = \left(-R(-\pi/2) + R(-\pi/2+\tau) \right)$$
 \end{thm}

\item \underline{The pseudo $2D$ case}:

The Vlasov-Poisson system considered in this case is the following $2D$ system:
 
  \begin{eqnarray}
\label{pseudo2D-1}
 \partial_t f_\epsilon +  \frac{v}{\epsilon} .\nabla_x f_\epsilon + \left( E_\epsilon + \frac{v^{\perp}}{\epsilon} \right).\nabla_v f_\epsilon =0 \\
  f_{\epsilon {\vert t=0}} = f_0 \\
  E_\epsilon =  - \nabla V_\epsilon, -\Delta_x V_\epsilon = \rho_\epsilon \\
\label{pseudo2D-2}
  \rho_\epsilon = \int f_\epsilon dv
 \end{eqnarray}
 
 If $v=(v_x,v_y)$, $v^\perp$ is defined by $(v_y,-v_x)$.

We recall that there exist global weak solutions of Vlasov-Poisson systems in the sense of Arsenev (\cite{Ar}).

Assuming here that $f_0 \geq 0, f_0 \in L^1_{x,v} \cap L^p_{x,v}$ (for some $p>2$) and that the initial energy is bounded, Frénod and Sonnendrücker proved the following theorem (we voluntarily write an unprecise meta-version of the result)
 \begin{thm} 
 For each $\epsilon$, let $(f_\epsilon,E_\epsilon)$ be a solution in the sense of Arsenev to (\ref{pseudo2D-1})-(\ref{pseudo2D-2}).
 
Then, up to a subsequence, $f_\epsilon$ weakly converges to a function $f$
 Moreover, there exists a function $G$ such that :
 \begin{equation}
f=\int_0^{2\pi}G(t,x+\mathcal{R}(\tau)v,R(\tau)v)d\tau
 \end{equation}
 and $G$ satisfies :
  \begin{eqnarray*}
  \partial_t G + \frac{1}{2\pi}\left(\int_0^{2\pi} \mathcal{R}(\tau)\mathcal{E}(t,\tau,x+\mathcal{R}(-\tau)v)d\tau \right).\nabla_x G \\
  +\frac{1}{2\pi}\left(\int_0^{2\pi} R(\tau)\mathcal{E}(t,\tau,x+\mathcal{R}(-\tau)v)d\tau \right).\nabla_v G =0
 \end{eqnarray*}
 $$  G_{\vert t=0} =  f_0$$
 $$ \mathcal{E}= -\nabla \Phi, \text{          } -\Delta \Phi = \int G(t,x+\mathcal{R}(\tau)v,R(\tau)v)dv$$
denoting by $R$ and $\mathcal{R}$ the linear operators defined by :
  $$R(\tau)= \begin{bmatrix}\cos \tau & -\sin \tau  \\ \sin \tau & \cos \tau \end{bmatrix}, \mathcal{R}(\tau) = \left(R(-\pi/2) - R(-\pi/2+\tau) \right)$$
 \end{thm}
In this case, we have to introduce an additional variable, the ``fast-time'' variable $\tau$  which comes from the fact that we need to precisely describe the oscillations in order to study the limit in non-linear terms.

  \end{enumerate}
 Note that the authors actually developped a generic framework that allows them to deal with different scalings and to give a precise approximation at any order. We do not wish to do so in our study.
 
\subsection{A bit of homogenization theory and some useful definitions}

Let us now precisely state the ``2-scale'' convergence tools used in this paper.

 \begin{deft}
 Let $X$ be a separable Banach space, $X'$ be its topological dual space and $(.,.)$ the duality bracket between $X'$ and $X$.
 For all $\alpha>0$, denote by $\mathcal{C}_{\alpha}(\mathbb{R},X)$ (respectively $L^{q'}_{\alpha}(\mathbb{R} ;X')$) the space of $\alpha$-periodic continuous (respectively $L^{q'}$) functions on $\mathbb{R}$ with values in $X$.
 Let $q\in[1;\infty[$. 
 
 Given a sequence $(u_\epsilon)$ of functions belonging to the space $L^{q'}(0,t ;X')$ and a function $U^0(t,\theta) \in L^{q'}(0,T ; L^{q'}_{\alpha}(\mathbb{R} ;X'))$ we say that
 $$u_\epsilon \text{   2-scale converges to   } U^0$$
 if for any function $\Psi \in L^q(0,T ; \mathcal{C}_{\alpha}(\mathbb{R},X))$ we have:
 \begin{equation}
 \lim_{\epsilon \rightarrow 0} \int_0^T \left(u_\epsilon(t), \Psi\left(t,\frac{t}{\epsilon}\right) \right)dt=\frac{1}{\alpha} \int_0^T \int_0^{\alpha} \left(U^0(t,\tau),\Psi(t,\tau) \right) d\tau dt
 \end{equation}
 \end{deft}
 
 \begin{thm}\label{2scale}
 Given a sequence $(u_\epsilon)$ bounded in $L^{q'}(0,t ;X')$, there exists for all $\alpha>0$ a function $U^0_\alpha \in L^{q'}(0,T ; L^{q'}_{\alpha}(\mathbb{R} ;X'))$ such that up to a subsequence, 
  $$u_\epsilon \text{   2-scale converges to   } U^0_\alpha$$
  The profile $U^0_\alpha$ is called the $\alpha$-periodic two scale limit of $u_\epsilon$ and the link between $U^0_\alpha$ and the weak-* limit $u$ of $u_\epsilon$ is given by:
  \begin{equation}
  \frac{1}{\alpha}\int_0^\alpha U^0d\tau = u
  \end{equation}
 \end{thm}

We also introduce some notations:
\begin{notat} 

We define for all $p\in [1;\infty]$ the space $L^p_{x,v} {:=} L^p_x(\mathbb{T}^d,(L^p_v(\mathbb{R}^d)))$.

In the same fashion, we define the spaces $L^p_{t,x}$, $L^p_{t,x,v}$...

Let $L^p_{2\pi, \tau}$ be the space of $2\pi$-periodic functions of $\tau$ which are in $L^p_\tau$.

Let $L^p_{x,\text{loc}}$ be the space of functions $f$ such that for all infinitely differentiable cut-off functions $\varphi \in \mathcal{C}_c^\infty$, $\varphi f$ belongs to $L^p_x$. We will say that a sequence $(f_\epsilon)$ is uniformly bounded in $L^p_{x,\text{loc}}$ if for each compact set $K$, the sequence of the restrictions to $K$ is uniformly bounded in $L^p_x$ with respect to $\epsilon$ (but this bound can depend on $K$). 

We will also use the same notations for Sobolev spaces $W^{s,p}$ ($s \in \mathbb{R}$).
\end{notat}

 \subsection{Statement of the result}

 In this paper we prove that the $2$-scale convergence established in the previous $2D$ case is also true in our $3D$ framework. The difficulty comes from the fact that there is no uniform elliptic regularity for the electric field because of the factor $\epsilon^2$ in front of $\Delta_{x_\parallel}$ in the Poisson equation:
$$V_\epsilon- \epsilon^2 \Delta_{x_\parallel} V_\epsilon - \Delta_{x_\perp} V_\epsilon = \int f_\epsilon dv - \int f_\epsilon dv dx$$
 In particular there is no a priori regularity on $x_\parallel$ and therefore no strong compactness. Nevertheless, we actually prove that due to the particular form of the asymptotic equation, the moments of the solution  with respect to $v_\parallel$ are more regular in $x_\parallel$ than the solution itself.  We can then easily pass to the weak limit.

The reason why we have opted for this strange Poisson equation instead of the usual one will appear at the end of the next section and especially in the last one. Roughly speaking it allows us to ``kill'' plasma waves which appear in the parallel direction due to the quasineutrality.

 Notice that this result is in the same spirit as the proof of the weak stability of the Vlasov-Maxwell system by DiPerna and Lions (\cite{DPL}), where the authors have regularity on moments, by opposition to the proof of the weak stability of the Vlasov-Poisson system by Arsenev (\cite{Ar}), where the author has compactness on the electric field. Actually our result is a kind of a hybrid one, since we get on one hand regularity with respect to $x_\perp$ by elliptic regularity and in the  other hand regularity with respect to $x_\parallel$ by averaging.
 
  We assume here that the initial data $(f_{\epsilon,0})_{\epsilon >0}$ satisfy the following conditions: 
\begin{itemize}
 \item $f_{\epsilon,0}\geq 0$ (positivity)
 \item $(f_{\epsilon,0})_{\epsilon>0}$ is uniformly bounded with respect to $\epsilon$ in $L^1_{x,v} \cap L^p_{x,v}$ (for some $p>3$) and for each $\epsilon$, $\int f_{\epsilon, 0} dxdv =1$.
 \item The initial energy is uniformly bounded with respect to $\epsilon$:
$$\left(\int f_{\epsilon, 0} \vert v \vert^2 dvdx +  \epsilon \int V_{\epsilon, 0}^2 dx  + \epsilon\int \vert \nabla_{x_\perp} V_{\epsilon, 0} \vert^2 dx + \epsilon^3 \int \vert \nabla_{x_\parallel} V_{\epsilon, 0}\vert^2 dx\right) \leq C$$
\end{itemize}

\begin{thm}\label{principal}
  For each $\epsilon$, let $(f_\epsilon,E_\epsilon)$ in $L^\infty_t(L^1_{x,v} \cap L^p_{x,v}) \times L^\infty_t(L^2_x)$ be a global weak solution in the sense of Arsenev to (\ref{real3D}).
Then up to a subsequence we have the following convergence as $\epsilon$ tends to $0$:
 
 \begin{eqnarray}
 f_{\epsilon,0} &\text{weakly-* converges to }& f_0 \in L^p_{x,v} \\
 f_\epsilon &\text{2-scale converges to }& F \in L^\infty_t(L^\infty_{2\pi,\tau}(L^1_{x,v} \cap L^p_{x,v})) \\
 E_\epsilon &\text{2-scale converges to }& \mathcal{E} \in L^\infty_t(L^\infty_{2\pi,\tau}(L^{3/2}_{x_\parallel}(W^{1,\frac{3}{2}}_{x_\perp})))
 \end{eqnarray}
 Moreover, there exists a function $G \in L^\infty_t(L^1_{x,v} \cap L^p_{x,v})$ such that:
 \begin{equation}
 F(t,\tau,x,v)=G(t,x+\mathcal{R}(\tau)v,R(\tau)v)
 \end{equation}
 and $(G,\mathcal{E})$ is solution to:
  \begin{eqnarray*}
  \partial_t G + v_\parallel.\nabla_x G + \frac{1}{2\pi}\left(\int_0^{2\pi} \mathcal{R}(\tau)\mathcal{E}(t,\tau,x+\mathcal{R}(-\tau)v)d\tau \right).\nabla_x G \\
  +\frac{1}{2\pi}\left(\int_0^{2\pi} R(\tau)\mathcal{E}(t,\tau,x+\mathcal{R}(-\tau)v)d\tau \right).\nabla_v G =0
 \end{eqnarray*}
 $$  G_{\vert t=0} =   f_0$$
 $$ \mathcal{E}= (-\nabla_\perp  V,0), \text{          }V -\Delta_\perp V = \int G(t,x+\mathcal{R}(\tau)v,R(\tau)v)dv- \int f_0 dvdx$$
  denoting by $R$ and $\mathcal{R}$ the linear operators defined by:
 $$R(\tau)= \begin{bmatrix}\cos \tau & -\sin \tau & 0 \\ \sin \tau & \cos \tau & 0 \\ 0 & 0 & 1\end{bmatrix}, \mathcal{R}(\tau) = \left(R(-\pi/2) - R(-\pi/2+\tau) \right)$$

 \end{thm}
 
  As it has been said, for the proof of this theorem, we will first prove a proposition which gives the regularity of moments in $v_\parallel$ of the solution. For this, we use an averaging lemma. The beginning of the proof is very similar to the proof in the $2D$ case, but we will give it again for the sake of completeness.
 \par
\begin{rques}
\begin{enumerate}
 \item The assumption on the initial energy may, at first sight, look a bit restrictive but in the ``usual'' Vlasov-Poisson scaling, it only means that the inital electric potential and field are bounded in $L^2$.

\item The constant $q=3$ will come quite naturally from Lemma \ref{inter} and Proposition \ref{ave}.

\item This theorem implies that for a given non-negative initial data $G_{\vert t=0}=G_0$ in $L^1_{x,v}\cap L^p_{x,v}$ (with $p>3$) and satisfying the energy bound, the asymptotic system admits at least one global weak solution $G \in L^\infty_t(L^1_{x,v} \cap L^p_{x,v})$. With the additional assumptions on the inital data:
\begin{eqnarray*}
 G_0 \in W^{1,1}_{x,v},\\
 \Vert (1  + \vert v \vert^4) G_0 \Vert_{L^\infty_{x,v}} < \infty, \\
  \Vert (1  + \vert v \vert^4) DG_0 \Vert_{L^\infty_{x,v}} < \infty.
\end{eqnarray*}
we are actually able to prove the uniqueness of the solution, using the same ideas than Degond in \cite{De} (and also used afterwards by Saint-Raymond in a gyrokinetic context (\cite{SR1})). Hence, it means that if the whole sequence $(f_{\epsilon,0})$ satisfies the same additional estimates, uniformly with respect to $\epsilon$, and weak * converges to some $f_0$ then there is also convergence for the whole sequence $(f_\epsilon)$.

\end{enumerate}
\end{rques}

 \section{A priori uniform estimates for the scaled Vlasov-Poisson system}
 
\subsection{Conservation of $L^p$ norms and energy for the scaled system}

 In this section we give a priori estimates which are very classical for the Vlasov-Poisson system (used for example in \cite{FS1}, \cite{FS2}, \cite {GSR1}). In order to recall how one can get them, we will give some formal computations. If one wants to have rigorous proofs, one should deal with smooth and compactly supported functions, namely with a sequence $(f^n_\epsilon)_{n\geq0}$ of solutions of some regularized Vlasov-Poisson equations then pass to the limit (that is the way one can clasically build a global weak solution in the sense of Arsenev (\cite{Ar})).

First, as usual for such Vlasov equations, $L^p$ norms are conserved (we work here at a fixed $\epsilon$):
 
 \begin{lem}
 For all $1\leq p\leq \infty$, 
 \begin{equation}
 \forall t \geq 0, \Vert f(t)\Vert_{L^p_{x,v}} \leq \Vert f(0)\Vert_{L^p_{x,v}}
 \end{equation}
 Moreover, $f_0\geq 0$ if and only if $\forall t \geq 0, f(t)\geq 0$ (referred to as the maximum principle)
 \end{lem}
 That precisely means that if $f_0 \in L^p_{x,v}$, then $f \in L^\infty_t(L^p_{x,v})$.
 
 Let us now compute the energy for the scaled system:
 
 \begin{lem}
\label{NRJ}
 We have the estimate:
 \begin{equation}
\mathcal{E}_\epsilon(t)= \left(\int f_\epsilon \vert v \vert^2 dvdx +  \epsilon \int V_\epsilon^2 dx  + \epsilon\int \vert \nabla_{x_\perp} V_\epsilon \vert^2 dx + \epsilon^3 \int \vert \nabla_{x_\parallel} V_\epsilon\vert^2 dx\right)\leq \mathcal{E}_\epsilon(0)
  \end{equation}
  In particular if there exists $C>0$ independent of $\epsilon$ such that $\mathcal{E}_\epsilon(0)\leq C$, then:
  \begin{equation}
  \int f_\epsilon \vert v \vert^2 dvdx \leq C
  \end{equation}
 \end{lem}
 
 \begin{proof}[Formal proof]
We multiply the scaled Vlasov equation by $\vert v \vert^2$ and integrate with respect to $x$ and $v$.
\begin{eqnarray*}
 \int \partial_t f_\epsilon \vert v \vert^2 dv dx+ \int E_\epsilon.\nabla_v f_\epsilon \vert v \vert^2 dvdx &=& \frac{d}{dt}\left( \int f_\epsilon \vert v \vert^2 dv dx\right) - 2\int E_\epsilon(x).v f_\epsilon dvdx =0
\end{eqnarray*}

We then integrate the Vlasov equation  with respect to $v$. We get the so called conservation of charge:
\begin{equation}
 \frac{d}{dt}\left( \int fdv \right) + \nabla_{x_\parallel}.\left( \int fv_\parallel dv\right) + \frac{\nabla_{x_\perp}}{\epsilon}.\left( \int f v_\perp dv\right) = 0  
\end{equation}

Therefore, we have:
\begin{eqnarray*}
 \int E_\epsilon(x).v f_\epsilon dvdx &=& - \int (\nabla_{x_\perp} V_\epsilon, \epsilon \nabla_{x_\parallel} V_\epsilon).v f_\epsilon  dv dx
\\ &=&  \int V_\epsilon \left( \nabla_{x_\perp}.(f_\epsilon v_\perp) + \epsilon \nabla_{x_\parallel}.(f_\epsilon v_\parallel)\right) dv dx
\\ &=& - \epsilon \int V_\epsilon \partial_t f_\epsilon dv dx
\end{eqnarray*}

Finally, using the Poisson equation, we get:
\begin{eqnarray*}
 - \epsilon \int V_\epsilon \partial_t f_\epsilon dv dx &=& -\epsilon \int V_\epsilon \partial_t\left( V_\epsilon-\epsilon^2 \Delta_{x_\parallel} V_\epsilon -  \Delta_{x_\perp} V_\epsilon\right) dx
\\&=& -\epsilon \left(\int V_\epsilon \partial_t V_\epsilon dx+ \int \nabla_{x_\perp}V_\epsilon \partial_t\nabla_{x_\perp}V_\epsilon  dx + \epsilon^2 \int  \nabla_{x_\parallel}\partial_t V_\epsilon\nabla_{x_\parallel}V_\epsilon dx \right)
\\ &=&  -\epsilon \frac{1}{2}\frac{d}{dt}\left( \int V_\epsilon^2 dx + \int \vert \nabla_{x_\perp} V_\epsilon \vert^2 dx + \epsilon^2 \int \vert \nabla_{x_\parallel} V_\epsilon\vert^2 dx \right)
\end{eqnarray*}

Thus it comes:

 \begin{equation}
 \frac{d}{dt} \left(\int f_\epsilon \vert v \vert^2 dvdx  + \epsilon \int V_\epsilon^2 dx + \epsilon \int \vert \nabla_{x_\perp} V_\epsilon \vert^2 dx + \epsilon^3 \int \vert \nabla_{x_\parallel} V_\epsilon\vert^2 dx\right)=0
  \end{equation}

 \end{proof}

 \subsection{Regularity of the electric field}
 
 Let us recall a classical lemma obtained by a standard real interpolation argument:
 
  \begin{lem}
\label{inter}Let $f(x,v)$ be a mesurable positive function on $\mathbb{R}^3 \times \mathbb{R}^3$. Then:  
  \begin{equation}
  \int \left(\int f(x,v)dv\right)^{3/2}dx \leq C \Vert f\Vert_{L^3_{x,v}}^{3/4} \left( \int \vert v\vert^2 f dxdv\right)^{3/4}
  \end{equation}
 \end{lem}
 
\begin{proof}
 For any $R>0$, we can write the following decomposition:
\begin{eqnarray*}
 \int f_\epsilon dv &=& \int_{\vert v\vert \leq R} f dv + \int_{\vert v\vert >R} f dv \\
&\leq& CR^2\Vert f \Vert_{L^3_v} + \frac{1}{R^2}\int  \vert v\vert^2 f dv
\end{eqnarray*}
Then we can take $R$ such that $R^2\Vert f \Vert_{L^3_v}=\frac{1}{R^2}\int  \vert v\vert^2 f dv$ so that we get:
\begin{equation}
 \int f dv \leq C \left(\int f^3 dv \right)^{1/6}\left(\int \vert v \vert^2 f dv\right)^{1/2}
\end{equation}
We then raise the quantities to the power $3/2$, integrate with respect to $x$ and use Hölder's inequality which gives the estimate.

\end{proof}

 By conservation of the $L^3$ norm and the uniform bound on the intial energy,  Lemmas \ref{NRJ} and \ref{inter}  entail that:
 
 \begin{equation}
 \rho_\epsilon \in L^\infty_t(L^{3/2}_x)
 \end{equation}
 and the norm is bounded uniformly with respect to $\epsilon$.
 
 We now use the Poisson equation to compute the regularity of the electric field. Let us recall that:
 \begin{eqnarray*}
 E_\epsilon= \left(-\epsilon \nabla_{x_\parallel} V_\epsilon, - \nabla_{x_\perp} V_\epsilon \right) \\
V_\epsilon -\epsilon^2 \Delta_{x_\parallel} V_\epsilon -  \Delta_{x_\perp} V_\epsilon = \rho_\epsilon - \int \rho_0 dx
 \end{eqnarray*}
 
 \begin{lem}
\label{elec}
With the above notations and assumptions:

 $E_\epsilon$ is uniformly bounded with respect to $\epsilon$ in ${ L^\infty_t(L^{3/2}_{x_\parallel}(W^{1,3/2}_{x_\perp}))}$

 \end{lem}
 
\begin{proof} Let $\epsilon>0$ and $t>0$ be fixed. For the sake of simplicity we write $V$ instead of $V_\epsilon$ and $E$ instead of $E_\epsilon$.

For any function $f(x_\parallel,x_\perp)$, define the rescaled function $\tilde{f} (z,x_\perp)$ by 
$$\tilde{f} \left(\frac{x_\parallel}{\epsilon},x_\perp\right)=\epsilon^{\frac{2}{3}}f(x_\parallel,x_\perp)$$
so that:

\begin{equation}\label{ident}\Vert \tilde{f} \left(z,x_\perp\right)\Vert_{L^{3/2}_{z}} = \Vert f(x_\parallel,x_\perp) \Vert_{L^{3/2}_{x_\parallel}}
\end{equation}

The Poisson equation becomes:
$$\tilde{V} -\Delta_{z} \tilde{V} -  \Delta_{x_\perp} \tilde{V} =  \tilde{\rho}- \epsilon^{\frac{2}{3}}\int \rho_0 dx$$
and the scaled electric field is given by:
$$ \tilde{E}= \left(-\nabla_{z} \tilde{V} , - \nabla_{x_\perp} \tilde{V}  \right)$$

 Since $\rho_\epsilon(t,.,.)$ and $V_\epsilon$ are uniformly bounded in $L^{3/2}_x$, standard results of elliptic regularity on the torus $\mathbb{T}^2\times \frac{1}{\epsilon} \mathbb{T}$ show that there exists $C>0$ independent of $\epsilon$ such that:

$$\Vert \tilde{V} \Vert _{{W}^{2,3/2}_{z,x_\perp}}\leq C \left\Vert \tilde{\rho}- \epsilon^{\frac{2}{3}}\int \rho_0 dx\right\Vert_{L^{3/2}_{z,x_\perp}}$$

\begin{rque}
 Notice here that due to the dilatation of order $\frac{1}{\epsilon}$ in the parallel direction, being periodic in this direction does not make things easier.
\end{rque}

Thanks to (\ref{ident}) we get:

$$\Vert \tilde{V} \Vert _{{W}^{2,3/2}_{z,x_\perp}} \leq C \left\Vert \rho- \int \rho_0 dx\right\Vert_{L^\infty_t(L^{3/2}_x)}\leq C_0$$
with $C_0$ independent of $\epsilon$.

Consequently, we have:
$$
\Vert \tilde{E} \Vert_{L^{3/2}_z(W^{1,3/2}_{x_\perp})} \leq \Vert \tilde{E} \Vert_{W^{1,3/2}_{z,x_\perp}}\leq C_0
$$
Finally from (\ref{ident}) we get 
 $$\Vert E_\epsilon\Vert_{ L^\infty_t(L^{3/2}_{x_\parallel}(W^{1,3/2}_{x_\perp}))}\leq C_0$$

 \end{proof}
 
 We can see as expected that the regularity of the electric field with respect to the $x_\parallel$ variable is not sufficient to get some strong compactness. 
 
\begin{rques}
\begin{enumerate}
\item We can write the identity:
\begin{equation}
 -\Delta_{x_\perp} V_\epsilon=-\Delta_{x_\perp}(Id-\epsilon^2 \Delta_{x_\parallel}  -  \Delta_{x_\perp})^{-1}\left( \rho_\epsilon- \int \rho_\epsilon dx\right)
\end{equation}
so that, thanks to elliptic estimates on the torus $\mathbb{T}^2$, $V_\epsilon \in L^{3/2}_{x_\parallel}(W^{2,3/2}_{x_\perp})$. Consequently, $\partial_{x_\parallel} V_\epsilon$ is bounded in $L^{3/2}_{x_\perp}(W^{-1,3/2}_{x_\parallel})$. This implies that $E_{\epsilon,\parallel}=-\epsilon \partial_{x_\parallel} V_\epsilon$ tends to zero in the sense of distributions.

\item A typical function $\varphi_\epsilon$ such that $\varphi_\epsilon$ is bounded in $L^p$ and $\frac{1}{\epsilon}\varphi_\epsilon$ is bounded in $W^{-1,p}$ is the oscillating function $\cos(\frac{1}{\epsilon} x) $. This indicates that $E_{\epsilon,\parallel}$ oscillates with a frequency of order $\frac{1}{\epsilon}$ in the parallel direction.

  \item If we work with the usual Poisson equation $$ -\epsilon^2 \Delta_{x_\parallel} V_\epsilon -  \Delta_{x_\perp} V_\epsilon = \rho_\epsilon - \int \rho_\epsilon dx$$ we only get homogeneous estimates for $V_\epsilon$ and we have not been able to deal with such anisotropic estimates in the following of the paper (namely in the estimates of Proposition \ref{ave}).
Roughly speaking, if $V$ is a solution of the Poisson equation $-\Delta V= \rho$ with $\rho \in L^{3/2}(\mathbb{R}^3)$, we can only say that $V \in \dot{W}^{2,3/2}(\mathbb{R}^3)$ (the homogeneous Sobolev space) and not $W^{2,3/2}$. 

\item This difficulty seems to be not only a technical one, but appears to be linked to the existence of plasma waves (with frequence and magnitude of order $\frac{1}{\sqrt{\epsilon}}$) in the parallel direction which prevents us from passing directly to the limit $\epsilon\rightarrow 0$ (see \cite{Gre} and last section).

\end{enumerate}

\end{rques}

 \section{Proof of Theorem \ref{principal}}

\begin{proof}
The first two steps are identical to the one given in \cite{FS2}. For the sake of completeness we recall here the main arguments and refer to \cite{FS2} for the details.

\subsection*{\underline{Step 1}: Deriving the constraint equation}

First of all, since $(f_\epsilon)$ is bounded in $L^\infty_t(L^1_{x,v}\cap L^p_{x,v})$, Theorem \ref{2scale} shows that for all $\alpha>0$:
$$f_\epsilon \text{   2-scale converges to   } F_\alpha \in L^\infty(0,T ; L^\infty_\alpha (\mathbb{R} ; L^p_{x,v}))$$

Let $\Psi(t,\tau,x,v)$ be an $\alpha$-periodic oscillating test function in $\tau$ and define: $$\Psi^\epsilon \equiv \Psi(t, \frac{t}{\epsilon}, x,v)$$ 
We start by writing the weak formulation of the scaled Vlasov equation against $\Psi^\epsilon$. Since $$\nabla_{x_\parallel}. v_\parallel = \nabla_{x_\perp}. v_\perp = \operatorname{div}_v\left( E_\epsilon + \frac{v\wedge e_z}{\epsilon} \right)=0,$$ we get the following equation:

\begin{align*}
 \int f_\epsilon \left(  (\partial_t \Psi)^\epsilon + \frac{1}{\epsilon} (\partial_\tau \Psi)^\epsilon + v_\parallel .(\nabla_x \Psi)^\epsilon +  \frac{v_\perp}{\epsilon} .(\nabla_x \Psi)^\epsilon + \left( E_\epsilon + \frac{v\wedge e_z}{\epsilon} \right).(\nabla_v \Psi)^\epsilon \right)dt dx dv \\= -\int f_0 \Psi(0,0,x,v)dx dv 
\end{align*}

Multiply then by $\epsilon$ and pass up to a subsequence to the ($2$-scale) limit. We get the so called constraint equation for the $\alpha$-periodic profile $F_\alpha$:

\begin{equation}
\partial_\tau F_\alpha + v_\perp.\nabla_x F_\alpha + v\wedge e_z.\nabla_v F_\alpha=0,
\end{equation}
which means that $F_\alpha$ is constant along the characteristics:
\begin{eqnarray}
\frac{dV}{d\tau}&=&V\wedge e_z\\
\frac{dX}{d\tau}&=&V_\perp
\end{eqnarray}

A straightforward calculation therefore shows that there exists $F^0_\alpha \in L^\infty(0,T, L^p_{x,v})$ such that:
\begin{equation}
F_\alpha(t,\tau,x,v)= F^0_\alpha(t,x+\mathcal{R}(\tau)v,R(\tau)v)
\end{equation}
with:
 $$R(\tau)= \begin{bmatrix}\cos \tau & -\sin \tau & 0 \\ \sin \tau & \cos \tau & 0 \\ 0 & 0 & 1\end{bmatrix}, \mathcal{R}(\tau) = \left(-R(-\pi/2) + R(-\pi/2+\tau) \right)$$

i.e. $\mathcal{R}(\tau)=\begin{bmatrix}\sin \tau & \cos \tau - 1 & 0 \\ 1- \cos \tau & \sin \tau & 0 \\ 0 & 0 & 0\end{bmatrix}$.

Since $R$ and $\mathcal{R}$ are $2\pi$-periodic, we will consider the $2\pi$ profile: indeed if $\alpha$ and $2\pi$ were incommensurable, $F_\alpha$ could not depend on $\tau$ and consequently we would have no information on the oscillations.

\subsection*{\underline{Step 2}: Filtering the essential oscillation}

We now look for the equation satisfied by $F^0_{2\pi}:=G$; we introduce the filtered function $g_\epsilon$:

\begin{equation}
g_\epsilon(t,x,v) =f_\epsilon(t,x+\mathcal{R}(-t/\epsilon)v,R(-t/\epsilon)v)
\end{equation}
(meaning that we have removed the oscillations)

We easily compute the equation satisfied by $g_\epsilon$:
   \begin{eqnarray}
\label{main}
  \partial_t g_\epsilon + v_\parallel.\nabla_x g_\epsilon +  \mathcal{R}(t/\epsilon)E_\epsilon(t, x+\mathcal{R}(-t/\epsilon)v).\nabla_x g_\epsilon \\
  +R(t/\epsilon)E_\epsilon(t, x+\mathcal{R}(-t/\epsilon)v).\nabla_v g_\epsilon =0 \nonumber
 \end{eqnarray}

\begin{rque}
Note here that $g_\epsilon$ 2-scale converges to $G$, and since it does not depend on $\tau$, it also weakly converges to $G$.
\end{rque}

\subsection*{\underline{Step 3}: Getting some regularity on moments}
From now on, the goal is to get some compactness for the moments of $g_\epsilon$ with respect to $v_\parallel$. The main tool we have in mind is the following averaging lemma proved by Bézard in \cite{BZ}, which is a refined version of the fundamental result of DiPerna, Lions and Meyer (\cite{DPLM}):

\begin{thm}
\label{Lp'temporel}Let $1<p\leq 2$. Let $f,g \in
L^p(dt\otimes dx\otimes dv)$ be solutions of the following transport equation
\begin{equation}
\label{TL}
 \partial_t f + v.\nabla_x f  = (I-\Delta_{t,x})^{\tau/2}(I-\Delta_v)^{m/2}g
\end{equation}
with $m\in \mathbb{R}^+, \tau \in [0,1[$. Then $\forall \Psi \in
\mathcal{C}^{\infty}_c(\mathbb{R}^d)$, $\rho_\Psi(t,x)=\int
f(t,x,v)\Psi(v)dv \in {W}^{s,p}_{t,x}(\mathbb{R}\times\mathbb{R}^d)$ where
\begin{equation}
s=\frac{1-\tau}{(1+m)p'}
\end{equation}
 Moreover,
\begin{equation}\Vert \rho_\Psi \Vert_{W^{s,p}_{t,x}(\mathbb{R}\times\mathbb{R}^d)} \leq C\left( \Vert f \Vert_{L^p(dt\otimes dx\otimes dv)} +
\Vert  g \Vert_{L^p(dt\otimes dx\otimes dv)}\right)\end{equation} ($C$ is a positive constant independent of $f$ and $g$)
\end{thm}

Averaging lemmas are an important feature of transport equations: since the transport equation (\ref{TL}) is hyperbolic, one can obviously not expect the solution $f$ to be more regular than the right hand side or the inital data. Nevertheless, if one considers the averaged quantity $\rho_\Psi$, one can actually notice a gain of regularity. This phenomenon was first observed  independently by Golse, Perthame and Sentis (\cite{GPS}) and Agoshkov (\cite{Ago}) then was formulated in a precise way for the first time by Golse, Lions, Perthame and Sentis (see \cite{GLPS}); it is referred to as ``velocity averaging''. There exists many refined versions of these results and numerous interesting applications in kinetic theory, but we shall not dwell on that.
We simply point out that this tool has been successfully applied to Vlasov equations, for instance to prove the existence of global weak solutions to the Vlasov-Maxwell system as it has been done by DiPerna and Lions (\cite{DPL}).
\par

 These results have been proved for functions with values in $\mathbb{R}$. Here, for our purpose, we need a new version of $L^p$ averaging lemma for functions with values in some Sobolev space  $W^{\lambda,p}(\mathbb{R}^k)$ ($k \in \mathbb{N}^*$). We prove the following result, which is sufficient in our case (probably an analogous of Bezard's optimal result is also true):
 
 \begin{lem}
 \label{abstractLp}
 Let $1<p < +\infty$ and $\lambda \in \mathbb{R}$. Let $f,g \in
L^p_{t,x,v}(W^{\lambda,p}_y)$ be solutions of the following transport equation
\begin{equation}
 \partial_t f + v.\nabla_x f  = (I-\Delta_v)^{m/2}g
\end{equation}
with $m\in \mathbb{R}^+$. Then $\forall \Psi \in
\mathcal{C}^{\infty}_c(\mathbb{R}^d)$, $\rho_\Psi(t,x)=\int
f(t,x,v)\Psi(v)dv \in {W}^{s,p}_{t,x}(W^{\lambda,p}_y)$ for any $s$ such that
\begin{equation}
s\leq s_2=\frac{1}{2(1+m)} \text{   for   } p= 2
\end{equation}
and
\begin{equation}
s<s_p=\frac{1}{(1+m)p'} \text{   for   } p\neq 2
\end{equation}
 Moreover,
\begin{equation}\Vert \rho_\Psi \Vert_{W^{s,p}_{t,x}(W^{\lambda,p}_y)} \leq C\left( \Vert f \Vert_{L^p_{t,x,v}(W^{\lambda,p}_y)} +
\Vert  g \Vert_{L^p_{t,x,v}(W^{\lambda,p}_y)}\right)\end{equation} ($C$ is a positive constant independent of $f$ and $g$)
 \end{lem}
 
 \begin{proof}[Sketch of proof] We prove the result in the stationnary case only:
  \begin{equation}
 v.\nabla_x f  = (I-\Delta_v)^{m/2}g
\end{equation}
 By standard arguments (see \cite{GLPS}) the general case then follows.
 
 The following estimate is obvious for $q=1$ or $q=+\infty$ (and actually we can not expect any smoothing effect) :
 \begin{equation}\Vert \rho_\Psi \Vert_{L^{q}_{x}(W^{\lambda,q}_y)} \leq C\left( \Vert f \Vert_{L^q_{x,v}(W^{\lambda,q}_y)} +
\Vert  g \Vert_{L^q_{x,v}(W^{\lambda,q}_y)}\right)\end{equation} 

For $p=2$, we prove the result as in Golse-Lions-Perthame-Sentis \cite{GLPS}. 
We denote by $\xi$ (resp. $\eta$) the Fourier variable associated to $x$ (resp. $y$).

The only point is to notice (using Fubini's inequality):
\begin{eqnarray*}
\Vert \rho_\Psi \Vert_{H^{s}_{x}(H^{\lambda}_y)}^2 &=& \int (1+\vert\xi\vert^2)^{s/2} \int (1+\vert\eta\vert^2)^{\lambda/2} \left( \int \mathcal{F}_{\xi,\eta} f \Psi(v)dv \right)^2 d\eta d\xi
\end{eqnarray*}
The proof is then identical and we get for $s=\frac{1}{2(1+m)}$:
 \begin{equation}\Vert \rho_\Psi \Vert_{H^{s}_{x}(H^{\lambda}_y)} \leq C\left( \Vert f \Vert_{L^2_{x,v}(H^{\lambda}_y)} +
\Vert  g \Vert_{L^2_{x,v}(H^{\lambda}_y)}\right)\end{equation} 

Finally the general case $1<p<+\infty$ is obtained by complex interpolation \cite{BL}.

 \end{proof}
 
  Equipped with this tool, we can now prove that moments in $v_\parallel$ are more regular with respect to $t$ and $x_\parallel$ than the solution itself.
 
 \begin{prop}
\label{ave}For each $\epsilon >0$, let $g_\epsilon$ be a function in $L^1_{x,v}\cap L^p_{x,v}$ (with $p>3$) bounded uniformly with respect to $\epsilon$ and satisfying:
   \begin{eqnarray*}
  \partial_t g_\epsilon + v_\parallel.\nabla_x g_\epsilon +  \mathcal{R}(t/\epsilon)E_\epsilon(t, x+\mathcal{R}(-t/\epsilon)v).\nabla_x g_\epsilon \\
  +R(t/\epsilon)E_\epsilon(t, x+\mathcal{R}(-t/\epsilon)v).\nabla_v g_\epsilon =0
 \end{eqnarray*}
 with $E_\epsilon$ the electric field uniformly bounded in $L^\infty_t(L^{3/2}_{x})$.
 
 Let $\Psi \in \mathcal{D}(\mathbb{R})$.
 %\Phi \in \mathcal{D}(\mathbb{R}^2), \
 Define $$\eta_\epsilon (t,x,v_\perp)= \int g_\epsilon(t,x,v)\Psi(v_\parallel)dv_\parallel$$
 %[ \mathcal{R}(-t/\epsilon).\nabla_x \Phi] (x_\perp-\mathcal{R}(t/\epsilon)v)
 Then,
 \begin{equation}
 \eta_\epsilon \text{  is uniformly bounded in  } W^{s, \gamma}_{t,x_\parallel, \text{loc}}(W^{-1, \gamma}_{x_\perp,v_\perp,\text{loc}})
 \end{equation}
 for  $\gamma \in ]1;2[$ defined by $\frac{1}{\gamma}=\frac{2}{3}+\frac{1}{p}$  and some $s\in]0;1[$ (depending on $\gamma$)

 \end{prop}
 
 \begin{proof}
 \begin{itemize}
 
\item The first step is to localize the equation. Let $K$ be the cartesian product of compact sets:
$$K=[0,{T}]\times {K}_{x_\parallel}\times {K}_{x_\perp}\times {K}_{v_\parallel}\times {K}_{v_\perp}$$
%$$\tilde{K}=[0,\tilde{T}]\times \tilde{K}_{x_\parallel}\times \tilde{K}_{x_\perp}\times \tilde{K}_{v_\parallel}\times \tilde{K}_{v_\perp}$$
 We now consider some positive smooth function $\Phi(t,x_\parallel,x_\perp,v_\parallel,v_\perp)$ which is $\mathcal C^\infty_c$ and which satisfies the condition:
 \begin{eqnarray}
 \Phi \equiv 0 \text{  outside  } K 
  %\Phi \equiv 0 \text{  outside  } {K} 
 \end{eqnarray}
 
 Noticing that:
$$\operatorname{div}_x \left(\mathcal{R}(t/\epsilon)E_\epsilon(t, x+\mathcal{R}(-t/\epsilon)v)\right) + \operatorname{div}_v \left(R(t/\epsilon)E_\epsilon(t, x+\mathcal{R}(-t/\epsilon)v)\right) =0$$
The equation satisfied by $g_\epsilon \Phi$ is the following one:
   \begin{eqnarray*}
  \partial_t (g_\epsilon \Phi) + v_\parallel.\nabla_x (g_\epsilon \Phi) =  - \underbrace{\nabla_x.(\mathcal{R}(t/\epsilon)E_\epsilon(t, x+\mathcal{R}(-t/\epsilon)v) g_\epsilon \Phi)}_{(1)} \\
  -\underbrace{\nabla_v.(R(t/\epsilon)E_\epsilon(t, x+\mathcal{R}(-t/\epsilon)v) g_\epsilon\Phi)}_{(2)} - \partial_t (\Phi) g_\epsilon -  v_\parallel.\nabla_x (\Phi) g_\epsilon \\
  +  \mathcal{R}(t/\epsilon)E_\epsilon(t, x+\mathcal{R}(-t/\epsilon)v).\nabla_x (\Phi)g_\epsilon  + R(t/\epsilon)E_\epsilon(t, x+\mathcal{R}(-t/\epsilon)v).\nabla_v (\Phi)g_\epsilon 
 \end{eqnarray*}
 
The idea is now to consider this equation as a kinetic equation with respect to the variables $(t,x_\parallel,v_\parallel)$ and with values in an abstract Banach space (which will be $W^{-1,\gamma}_{x_\perp,v_\perp}$).
 We then only study the first two terms of the right-hand side (noticing that the other terms have more regularity than these ones).

From now on, for the sake of simplicity and readability, we will write $L^p$ and $W^{s,p}$ norms without always specifying that they are taken on the compact support of $\Phi$.

 \item \subsubsection*{Estimate on the first term (1)}
Since $E_\epsilon$ does not depend on v, we have:
 $$E_\epsilon \in L^\infty_t(L^{3/2}_{x_\parallel}(L^{\infty}_v(L^ {3/2}_{x_\perp})))$$
 In particular if we restrict to compact supports:
 $$E_\epsilon \in L^{3/2}_{t,x,v}$$
 The second point is that the differential operator applied in (1) involves only derivatives with respect to the $x_\perp$ variable and not in the parallel direction: this remark is fundamental for using the averaging lemma \ref{abstractLp} (indeed, the case of a full derivative in $x_\parallel$ can not be handled). 

Hölder's inequality simply implies that:
 
 \begin{equation}
 \Vert \mathcal{R}(t/\epsilon)E_\epsilon(t, x+\mathcal{R}(-t/\epsilon)v)g_\epsilon \Phi\Vert_{L^\gamma_{x_\perp, v_\perp}}
 \leq 
 \Vert  E_\epsilon(t, x+\mathcal{R}(-t/\epsilon)v) \Phi\Vert_{L^{3/2}_{x_\perp, v_\perp}} \Vert g_\epsilon \Vert_{L^p_{x_\perp, v_\perp}}
\end{equation}
where  $\frac{1}{\gamma}=\frac{2}{3}+\frac{1}{p}$. Hence:
 \begin{equation}
 \Vert \nabla_x.(\mathcal{R}(t/\epsilon)E_\epsilon(t, x+\mathcal{R}(-t/\epsilon)v)g_\epsilon \Phi)\Vert_{W^{-1,\gamma}_{x_\perp, v_\perp}}
 \leq 
 \Vert  E_\epsilon(t, x+\mathcal{R}(-t/\epsilon)v) \Phi\Vert_{L^{3/2}_{x_\perp, v_\perp}} \Vert g_\epsilon \Vert_{L^p_{x_\perp, v_\perp}}
\end{equation}

Notice that the change of variables $(x,v)\mapsto(x+\mathcal{R}(s)v,v)$ has unit Jacobian for all $s \in \mathbb{R}$, so that:

\begin{equation}
  \Vert  E_\epsilon(t, x+\mathcal{R}(-t/\epsilon)v) \Phi\Vert_{L^{3/2}_{x_\perp, v_\perp}}=\Vert E_\epsilon(t, x)\Vert_{L^{3/2}_{x_\perp,v_\perp}}
\end{equation}

So finally we have, after integrating in $t,x_\parallel,v_\parallel$ and thanks to Hölder's inequality:

 \begin{eqnarray*}
 \Vert \nabla_x.\left(\mathcal{R}(t/\epsilon)E_\epsilon(t, x+\mathcal{R}(-t/\epsilon)v) g_\epsilon \Phi\right)\Vert_{L^\gamma_{t,x_\parallel,v_\parallel}(W^{-1,\gamma}_{x_\perp,v_\perp} )}\\
 \leq C\Vert E_\epsilon(t, x)\Vert_{L^{3/2}_{t,x_\parallel,v_\parallel}(L^{3/2}_{x_\perp,v_\perp})}\Vert g_\epsilon \Vert_{L^{p}_{t,x_\parallel,v_\parallel}(L^p_{x_\perp,v_\perp})}
\end{eqnarray*}

and $C$ is a constant independent of $\epsilon$.

\begin{rque}The regularity of $(1)$ with respect to $v_\perp$ is not optimal (since it involves no derivative in $v_\perp$ for $g_\epsilon$). Nevertheless we are interested in the regularity of the whole right hand side, and we will see that the term $(2)$ has this regularity in $v_\perp$.\end{rque}

 \item \subsubsection*{Estimate on the second term (2)}
  
  By the same method one gets:
   \begin{eqnarray*}
\Vert \nabla_v.\left(R(t/\epsilon)E_\epsilon(t, x+\mathcal{R}(-t/\epsilon)v)g_\epsilon \Phi\right)\Vert_{L^\gamma_{t,x_\parallel}( W^{-1,\gamma}_{v_\parallel}(W^{-1,\gamma}_{x_\perp,v_\perp}))}
 \\
 \leq C\Vert E_\epsilon(t, x) \Vert_{L^{3/2}_{t,x_\parallel}(L^{3/2}_{x_\perp,v})}\Vert g_\epsilon \Phi\Vert_{L^p_{t,x_\parallel}(L^p_{x_\perp,v})}
\end{eqnarray*}
Finally we see that the right hand side is uniformly bounded in:
$$L^\gamma_{t,x_\parallel,\text{loc}}(W^{-1,\gamma}_{v_\parallel,\text{loc}}(W^{-1,\gamma}_{x_\perp,v_\perp,\text{loc}}))$$

 \item \subsubsection*{Regularity of the moments}
  By lemma \ref{abstractLp} , for all $\Psi\in \mathcal{C}^\infty_c$, the moment:
  
  $$\eta_\epsilon (t,x,v_\perp)= \int g_\epsilon(t,x,v)\Psi(v_\parallel)dv_\parallel$$
 %[ \mathcal{R}(-t/\epsilon).\nabla_x \Phi] (x_\perp-\mathcal{R}(t/\epsilon)v)
 is then uniformly bounded  in the space $W^{s,\gamma}_{t,x_\parallel,\text{loc}}(W^{-1,\gamma}_{x_\perp,v_\perp,\text{loc}})$ for any $s>0$ with $s<\frac{1}{2\gamma'}$.
    \end{itemize}
 \end{proof}
 
 We can now prove that the sequence of moments $\eta_\epsilon$ is compact in a space of distributions which is the dual of some space where the sequence $(E_\epsilon)$ is uniformly bounded.
 
 \begin{coro}
\label{coro}
 There exists $\theta\in ]0,1[$ and
 $\eta \in W^{s\theta, 3}_{t,x_\parallel, \text{loc}}(W^{-\theta, 3}_{x_\perp,v_\perp,\text{loc}})$ such that for all $\xi>0$, up to a subsequence:
 
 \begin{equation}
\eta_\epsilon  \rightarrow \eta \text{ strongly in } L^{ 3}_{t, \text{loc}}(L^{3}_{x_\parallel,\text{loc}}(W^{-\theta-\xi, 3}_{x_\perp,v_\perp,\text{loc}}))
\end{equation}

 \end{coro}
 
 \begin{proof}
 By assumption on the initial data, there exists $q>3$ such that $f_0 \in L^q_{x,v}$; thanks to the a priori $L^q$ estimate, we get $g_\epsilon \in L^\infty_t(L^q_{x,v})$. Define $\gamma$ by: $$\frac{1}{\gamma}=\frac{2}{3}+\frac{1}{q}$$  
 The previous lemma shows that for some $s>0$:
 $$\eta_\epsilon \in {W}^{s, \gamma}_{t,x_\parallel \text{loc}}(W^{-1, \gamma}_{x_\perp,v_\perp,\text{loc}}) \text{ uniformly in } \epsilon$$

 Since $g_\epsilon \in L^q_{t,\text{loc}}(L^q_{x,v})$ and $\Psi$ has compact support, we get by  Hölder's inequality:
 $$\eta_\epsilon \in L^q_{t,\text{loc}} (L^q_{x_\parallel}(L^q_{ x_\perp, v_\perp}))$$
Since $\frac{1}{\gamma}>\frac{2}{3}>\frac{1}{3}$ and $\frac{1}{q}<\frac{1}{3}$, there exists $\theta\in ]0,1[$ such that 
 $$\frac{1}{3} = \frac{1-\theta}{q} + \frac{\theta}{\gamma}$$
 By interpolation (\cite{BL}) we deduce that:
 $$\eta_\epsilon \in W^{s\theta, 3}_{t, x_\parallel \text{loc}}(W^{-\theta, 3}_{x_\perp,v_\perp,\text{loc}}))$$

This implies that:
$$\eta_\epsilon \in {W}^{s\theta, 3}_{t, \text{loc}}(L^3_{x_\parallel}(W^{-\theta, 3}_{x_\perp,v_\perp,\text{loc}})) \text{ uniformly in } \epsilon$$
$$\eta_\epsilon \in {L}^{3}_{t, \text{loc}}({W}^{s\theta, 3}_{x_\parallel, \text{loc}}(W^{-\theta, 3}_{x_\perp,v_\perp,\text{loc}})) \text{ uniformly in } \epsilon$$

 We then use the following refined interpolation result proved by Simon in \cite{SIM}, which is, roughly speaking, an anisotropic adaptation of the classical Riesz-Fréchet-Kolmogorov criterion for compactness in $L^p$:

 \begin{thm} Let $1\leq p \leq \infty$ and $s>0$. 
 Let $T>0$ and $X,B,Y$ be three Banach Spaces such that $X\subset B \subset Y$ and with $X$ compactly embedded in $B$. Let $F$ be a bounded set of $L^p_t([0,T], X)\cap W^{s,p}_t([0,T],Y)$. 
Then $F$ is relatively compact in $L^p_t([0,T],B)$.
 \end{thm}

  This entails, thanks to Sobolev's embeddings, that the sequence $(\eta_\epsilon)$ is strongly relatively compact in $L^{3}_{t, \text{loc}}(L^{3}_{x_\parallel,\text{loc}}(W^{-\theta-\xi, 3}_{x_\perp,v_\perp,\text{loc}}))$, for all $\xi >0$.

 \end{proof}

 From now on, we consider $\xi$ such that $\theta +\xi <1$, which is of course possible since $\theta <1$.
 
\par

 \begin{rque}
 Following the remark in Step 2 and by uniqueness of the limit in the sense of distributions, we get:
 $$\eta = \int G\Psi(v_\parallel)dv_\parallel$$
 \end{rque}

\subsection*{\underline{Step 4}: Passing to the weak limit}
 We will first need a technical lemma which is obtained directly from the 2-scale convergence of $E_\epsilon$.
 \begin{lem}
\label{lem2scale} Up to a subsequence,

\begin{itemize}
 \item $\mathcal{R}(t/\epsilon)E_\epsilon(t, x+\mathcal{R}(-t/\epsilon)v)$ 2-scale converges to 
$\frac{1}{2\pi}\int_0^{2\pi}\mathcal{R}(\tau )\mathcal{E}(t, \tau, x+\mathcal{R}(-\tau)v)d\tau \in L^\infty_t(L^\infty_{2\pi, \tau}(L^{3/2}_{x_\parallel}(W^{1,3/2}_{x_\perp}))$
\item $R(t/\epsilon)E_\epsilon(t, x+\mathcal{R}(-t/\epsilon)v)$  2-scale converges to  $\frac{1}{2\pi}\int_0^{2\pi} R(\tau )\mathcal{E}(t, \tau, x+\mathcal{R}(-\tau)v)d\tau \in L^\infty_t(L^\infty_{2\pi, \tau}(L^{3/2}_{x_\parallel}(W^{1,3/2}_{x_\perp}))$
\end{itemize}
 \end{lem}

\begin{proof}
 $E_\epsilon$ is uniformly bounded in $L^\infty_t(L^{3/2}_{x_\parallel}(W^{1,3/2}_{x_\perp}))$, so there exists $\mathcal{E} \in L^\infty_t(L^{3/2}_{x_\parallel}(W^{1,3/2}_{x_\perp})$ such that $E_\epsilon$ 2 scale converge to $\mathcal{E}$.

We take $\Psi(t,\tau,x)$ a $2\pi$-periodic w.r.t. $\tau$ test function and use the 2 scale convergence of $E_\epsilon$:

\begin{eqnarray*}
 \int \mathcal{R}(t/\epsilon)E_\epsilon(t, x+\mathcal{R}(-t/\epsilon)v). \Psi(t,t/\epsilon,x) dt dx \\ = \int E_\epsilon(t, x) .^t\mathcal{R}(t/\epsilon)\Psi(t,t/\epsilon,x-\mathcal{R}(-t/\epsilon)v) dt dx
 \end{eqnarray*}
\begin{eqnarray*}
&\rightarrow& \frac{1}{2\pi}\int \int_0^{2\pi}\mathcal{E}(t,\tau,x).^t\mathcal{R}(\tau)\Psi(t,\tau,x-\mathcal{R}(-\tau)v) dt d\tau dx \\
&=& \frac{1}{2\pi}\int \int_0^{2\pi}\mathcal{E}(t,\tau,x+\mathcal{R}(-\tau)v) .^t\mathcal{R}(\tau)\Psi(t,\tau,x) dt d\tau dx
 \end{eqnarray*}

The proof is the same for $R(t/\epsilon)E_\epsilon(t, x+\mathcal{R}(-t/\epsilon)v)$.

\end{proof}
 \par

Now, we  can write the weak formulation of the kinetic equation (\ref{main}) against a smooth test function of the form $\Phi(t,x,v_\perp)\Psi(v_\parallel)$ with compact support.
 If we can pass to the limit for such test functions, then by density it will be also the case for all test functions. 
 
 Noticing that $\operatorname{div}_{x} v_\parallel =0$ and that
$$\operatorname{div}_x \left(\mathcal{R}(t/\epsilon)E_\epsilon(t, x+\mathcal{R}(-t/\epsilon)v)\right) + \operatorname{div}_v \left(R(t/\epsilon)E_\epsilon(t, x+\mathcal{R}(-t/\epsilon)v)\right) =0$$

we get:
 
 \begin{eqnarray*}
  \int \Big(\partial_t (\Phi(t,x,v_\perp)\Psi(v_\parallel)) + v_\parallel.\nabla_x (\Phi\Psi) +  \mathcal{R}(t/\epsilon)E_\epsilon(t, x+\mathcal{R}(-t/\epsilon)v).\nabla_x (\Phi\Psi) \\
  +R(t/\epsilon)E_\epsilon(t, x+\mathcal{R}(-t/\epsilon)v).\nabla_v (\Phi\Psi) \Big)g_\epsilon dtdx_\perp dx_\parallel dv_\perp dv_\parallel \\
  = -\int u_0 \Phi(0,x,v_\perp)\Psi(v_\parallel) dx dv \\
 \end{eqnarray*}
 
 We can easily take weak limits in the linear part $\partial_t g_\epsilon + v_\parallel . \nabla_x g_\epsilon$.
 
 Consider now the ``non linear'' term:
 \begin{eqnarray*}
 \int \mathcal{R}(t/\epsilon)E_\epsilon(t, x+\mathcal{R}(-t/\epsilon)v).\nabla_x \Phi(t,x,v_\perp) g_\epsilon \Psi(v_\parallel)dtdx_\perp dx_\parallel dv_\perp dv_\parallel= \\
  \int \mathcal{R}(t/\epsilon)E_\epsilon(t, x+\mathcal{R}(-t/\epsilon)v).\nabla_x \Phi(t,x,v_\perp) \left(\int g_\epsilon \Psi(v_\parallel)dv_\parallel\right)dtdx_\perp dx_\parallel dv_\perp 
  \end{eqnarray*}
  
  The convergence of this term can be established by the strong/weak convergence principle. Nevertheless, we have to carefully use this technique to get the result and we will explicitly evaluate the difference: 
  \begin{align*}
\Big\vert   \int \mathcal{R}(t/\epsilon)E_\epsilon(t, x+\mathcal{R}(-t/\epsilon)v)&.\nabla_x \Phi(t,x,v_\perp) \eta_\epsilon dtdx_\perp dx_\parallel dv_\perp   \\
  -\int \frac{1}{2\pi}&\int_0^{2\pi}\mathcal{R}(\tau )\mathcal{E}(t, \tau, x+\mathcal{R}(-\tau)v)d\tau.\nabla_x \Phi \eta dtdxdv_\perp \Big\vert  \\
 \leq \Big\vert  \int \Big(\mathcal{R}(t/\epsilon)E_\epsilon(t, x+\mathcal{R}(-t/\epsilon)v)&- \frac{1}{2\pi}\int_0^{2\pi}\mathcal{R}(\tau )\mathcal{E}(t, \tau, x+\mathcal{R}(-\tau)v)d\tau \Big)  .\nabla_x \Phi \eta dtdx dv_\perp \Big\vert   \\
  +\Big\vert  \int \mathcal{R}(t/\epsilon)E_\epsilon(t, x+\mathcal{R}(-t/\epsilon)v)&.\nabla_x \Phi(t,x,v_\perp)\left( \eta_\epsilon-\eta\right) dtdx_\perp dx_\parallel dv_\perp \Big\vert
  \end{align*}

  The first term of the right hand side converges to zero because of the 2-scale convergence of $E_\epsilon$ (Lemma \ref{lem2scale}). We can control the second term by:
  \begin{equation}
  C\Vert E_\epsilon.\nabla_x \Phi   \Vert_{L^{3/2}_{t}(L^{3/2}_{x_\parallel}(W^{\theta+\xi,3/2}_{x_\perp,v_\perp}))} \Vert \eta_\epsilon -\eta  \Vert_{L^{ 3}_{t}( L^{3}_{x_\parallel}(W^{-\theta-\xi, 3}_{x_\perp,v_\perp}))}
  \end{equation}
  (these norms are actually taken on the compact support of $\Phi$ but we do not write it for the sake of simplicity)

  Using the fact that $E_\epsilon$ is uniformly bounded in $L^{3/2}_{t,\text{loc}}(L^{3/2}_{x_\parallel}(W^{\theta+\xi,3/2}_{x_\perp}))$ (this is an easy consequence of Lemma \ref{elec}) and Corollary \ref{coro},
  $$\Vert \eta_\epsilon -\eta  \Vert_{L^{ 3}_{t}([0,T], L^{3}_{x_\parallel}(K_{x_\parallel}, W^{-\theta-\xi, 3}_{x_\perp,v_\perp}(K_{x_\perp}\times K_{v_\perp})))}\rightarrow 0,$$
we can deduce that 
  $$\left\vert  \int \mathcal{R}(t/\epsilon)E_\epsilon(t, x+\mathcal{R}(-t/\epsilon)v).\nabla_x \Phi(t,x,v_\perp)\left( \eta_\epsilon-\eta\right) dtdx_\perp dx_\parallel dv_\perp \right\vert \rightarrow 0$$
  
The proof is of course the same for the other non-linear term: $$R(t/\epsilon)E_\epsilon(t, x+\mathcal{R}(-t/\epsilon)v).\nabla_v g_\epsilon$$

To conclude let us compute the asymptotic equation satisfied by the 2-scale limit of $V_\epsilon$ denoted by $V$. We take $\Psi(t,\tau,x)$ a $2\pi$-periodic w.r.t. $\tau$ test function. We write the weak formulation of the Poisson equation:

\begin{eqnarray*}
\int V_\epsilon  \nabla_{x_\parallel}\Psi(t,t/\epsilon,x)dtdx + \\
\epsilon^2 \int \nabla_{x_\parallel}V_\epsilon  \nabla_{x_\parallel}\Psi(t,t/\epsilon,x)dtdx +  \int \nabla_{x_\perp}V_\epsilon  \nabla_{x_\perp}\Psi(t,t/\epsilon,x)dtdx \\ =\int f_\epsilon(t,x,v)\Psi(t,t/\epsilon,x) dt dv dx - \int \left( \int f_0 dvdx\right)\Psi(t,t/\epsilon,x) dt dv dx 
\end{eqnarray*}

We then pass to the 2 scale limit:

\begin{eqnarray*}
 & & \frac{1}{2\pi}\int \int_0^{2\pi} V(t,\tau,x)  \nabla_{x_\perp}\Psi(t,\tau,x)d\tau dtdx +0 \\
&+& \frac{1}{2\pi}\int \int_0^{2\pi} \nabla_{x_\perp}V(t,\tau,x)  \nabla_{x_\perp}\Psi(t,\tau,x)d\tau dtdx = \frac{1}{2\pi}\int  \int_0^{2\pi} F(t,\tau, x,v)\Psi(t,\tau,x) dt dv dx 
\\&-& \frac{1}{2\pi}\int  \int_0^{2\pi}\left(\int f_0 dvdx\right)\Psi(t,\tau,x) d\tau dv dx \\
 &=& \frac{1}{2\pi}\int  \int_0^{2\pi}G(t,\tau, x+\mathcal{R}(\tau)v, R(\tau)v)\Psi(t,\tau,x) d\tau dv dx \\&-& \frac{1}{2\pi}\int  \int_0^{2\pi}\left(\int f_0 dvdx\right)\Psi(t,\tau,x) d\tau dv dx
\end{eqnarray*}
from which we get the ``Poisson'' equation given in Theorem \ref{principal}:
$$V -\Delta_\perp V = \int G(t,x+\mathcal{R}(\tau)v,R(\tau)v)dv - \int f_0 dvdx$$

Moreover since $E_{\epsilon,\perp}=-\nabla_{x_\perp} V_\epsilon$ and thanks to Remark 1 following Lemma \ref{elec}, we easily get if we pass to the two-scale limit:
$$\mathcal{E}=(-\nabla_{x_\perp} V,0)$$

  \end{proof}
\section{Concluding comments}

\subsection{Comments on the result}

Finally we can see as in \cite{FS2} (namely by performing the change of variables $x=x_c-v^\perp$ and looking at the new equations in the so-called gyro-variables $(t,x_c,v)$) that the drift involving the electric field in the asymptotic ``kinetic'' equation corresponds to the electric drift that we mentioned in the introduction and which was expected to appear.
We also notice that the Poisson equation we get in the asymptotic system is the same than the one used in the numerical simulations of tokamak plasmas (see for example the GYSELA code in \cite{Gra}). Nevertheless, physicists do not get it in the same formal way: they claim that it only expresses the quasineutrality of the plasma (there is no ``real'' Poisson equation involved) and the perpendicular laplacian happens to appear due from the so-called ``polarization drift'' (\cite{Gra}, see also \cite{Wes} for a physical reference on the subject). It would be interesting to justify such a computation from a mathematical point of view.

At last, we wish to point out a really unpleasant feature of our model, which is that there is no parallel dynamics.

\subsection{An alternative model}

Let us give some comments on the gyrokinetic approximation of the system (\ref{fix}) which consists in considering a population of electrons in a fixed background of ions:
\begin{equation}
 n^i_\epsilon= \int f_0 dxdv
\end{equation}

Actually, quite surprisingly, this model engenders more interesting physical properties: in this case the parallel component of the electric field does not vanish but appears as a pressure in the end (which may bring difficulties both in the study of the asymptotic system and in numerical simulations).

  \begin{equation}
\label{fix}
\left\{
    \begin{array}{ll}
    \partial_{t} f_\epsilon + \frac{v_\perp}{\epsilon}.\nabla_{x} f_\epsilon + v_\parallel.\nabla_{x} f_\epsilon + (E_\epsilon+ \frac{v\wedge B}{\epsilon}).\nabla_{v} f_\epsilon = 0 \\
  E_\epsilon= (-\nabla_{x_\perp} {V}_\epsilon, -\epsilon\nabla_{x_\parallel} {V}_\epsilon)  \\
  -\epsilon^2\Delta_{x_\parallel} {V}_\epsilon -\Delta_{x_\perp} {V}_\epsilon = \int f_\epsilon dv - \int f_0 dvdx\\
     f_{\epsilon, t=0}=f_{\epsilon,0}
\end{array}
  \right.
\end{equation}

With the same computations as the present paper, we get:
 \begin{eqnarray}
 f_\epsilon &\text{2-scale converges to }& F \\
 E_\epsilon &\text{2-scale converges to }& \mathcal{E} 
 \end{eqnarray}
 In a formal sense, there exists a function $G$ such that:
 \begin{equation}
 F(t,\tau,x,v)=G(t,x+\mathcal{R}(\tau)v,R(\tau)v)
 \end{equation}
 and $(G,\mathcal{E})$ is solution to:
  \begin{eqnarray*}
  \partial_t G + v_\parallel.\nabla_x G + \frac{1}{2\pi}\left(\int_0^{2\pi} \mathcal{R}(\tau)\mathcal{E}(t,\tau,x+\mathcal{R}(-\tau)v)d\tau \right).\nabla_x G \\
  +\frac{1}{2\pi}\left(\int_0^{2\pi} R(\tau)\mathcal{E}(t,\tau,x+\mathcal{R}(-\tau)v)d\tau \right).\nabla_v G =0
 \end{eqnarray*}
 $$  G_{\vert t=0} =   f_0$$
 $$ \mathcal{E}= (-\nabla_\perp  V,\mathcal{E}_\parallel), \text{          } -\Delta_\perp V = \int G(t,x+\mathcal{R}(\tau)v,R(\tau)v)dv- \int f_0 dvdx$$
 still denoting by $R$ and $\mathcal{R}$ the linear operators defined by:
 $$R(\tau)= \begin{bmatrix}\cos \tau & -\sin \tau & 0 \\ \sin \tau & \cos \tau & 0 \\ 0 & 0 & 1\end{bmatrix}, \mathcal{R}(\tau) = \left(R(-\pi/2) - R(-\pi/2+\tau) \right)$$

The parallel component $\mathcal{E}_\parallel$ has to be seen as a pressure (or the Lagrange multiplier) associated to the ``incompressibility'' constraint $\int_{\mathbb{T}^2} \int G(t,x+\mathcal{R}(\tau)v,R(\tau)v)dvdx_\perp = \int f_0 dvdx$

Let us just give a few words on the difficulties that arise with this model. The Poisson equation can be restated as:
\begin{equation}
 -\epsilon^2\Delta_{x_\parallel} {V}_\epsilon -\Delta_{x_\perp} {V}_\epsilon = \int f_\epsilon dv - \int f_\epsilon dvdx_\perp + \int f_\epsilon dvdx_\perp - \int f_0 dvdx
\end{equation}
so that thanks to the linearity of the Poisson equation we can study separately two equations. The first one states:
\begin{equation}
 -\epsilon^2\Delta_{x_\parallel} {V}^1_\epsilon -\Delta_{x_\perp} {V}^1_\epsilon = \int f_\epsilon dv - \int f_\epsilon dvdx_\perp 
\end{equation}
For this part of the electric potential we get the same estimates as in lemma \ref{elec}. Indeed, $\int\left(\int f_\epsilon dv - \int f_\epsilon dvdx_\perp \right) dx_\perp=0$ so that we can use elliptic estimates on the torus $\mathbb{T}^2_{x_\perp}$. Consequently this electric potential does not give birth to any parallel dynamics, like in Theorem \ref{principal}.

The second one is:
\begin{equation}
 -\epsilon^2\Delta_{x_\parallel} {V}^2_\epsilon  = \underbrace{\int f_\epsilon dvdx_\perp - \int f_0 dvdx}_{\text{only depends on } x_\parallel}
\end{equation}

This equation, associated to the one giving the electric field $E^2_{\epsilon, \parallel}=-\epsilon \partial_{x_\parallel} V^2_\parallel$, is similar to the one studied by Grenier in \cite{Gre}, coupled to a Vlasov equation describing a quasineutral plasma. In this case it was shown that there exist plasma waves with both temporal and spatial oscillations with frequency $\frac{1}{\sqrt{\epsilon}}$ and magnitude of order $\frac{1}{\sqrt{\epsilon}}$. Because of these waves, it is much more difficult to pass to the limit in order to get a kinetic equation. Grenier managed to prove the convergence only for distribution functions with special form and got in the end a Euler-like system with an electric field interpreted as a Lagrange multiplier. Hence, in our case, this part of the electric field may engender an non-zero $\mathcal{E}_{\parallel}$.

For these reasons, it seems much harder to expect to prove a rigorous result similar to Theorem \ref{principal} with such a model.

\subsection{Prospects}

A way to pass to the limit in this latest case would be to use a relative entropy method like in the papers of Brenier (\cite{Br}) and Golse and Saint-Raymond (\cite{GSR2}). This will be the object of a forthcoming work. 

Another interesting issue would be to consider a ``true'' Boltzmann-Maxwell distribution for the electrons (not linearized like in this paper) and perform an asymptotic analysis, maybe also with a relative entropy method.

\begin{ack} The author would like to thank L. Saint-Raymond for having given him the idea to write this paper and for her careful reading and advice. The author is also indebted to M. Hauray for having pointed out to him the existence of plasma waves and for indicating Grenier's paper.
\end{ack}

\bibliographystyle{plain}
\bibliography{Larmor3D}

\end{document}